\documentclass[12pt]{amsart}

\usepackage{amsmath}
\usepackage{amsfonts}
\usepackage{graphicx}
\usepackage{latexsym}
\usepackage{amsbsy}  

\newtheorem{proposition}{Proposition}

\newtheorem{theorem}[proposition]{Theorem}

\newtheorem{definition}[proposition]{Definition}

\newcommand{\geoG}{\boldsymbol G}
\newcommand{\geoT}{\boldsymbol T}
\newcommand{\repG}{\sf G}
\newcommand{\repT}{\sf T}

\newcommand{\lie}{\mbox{Lie}}

\newcommand{\cplx}{\mathbb{C}}
\newcommand{\real}{\mathbb{R}}
\newcommand{\fps}{\mathbb{C}[[t]]}
\newcommand{\fls}{\mathbb{C}(\!(t)\!)}
\newcommand{\ggeq}{\succeq}
\newcommand{\gleq}{\preceq}
\newcommand{\glt}{\prec}
\newcommand{\ddim}{\mbox{dim}}
\newcommand{\Irr}{\mbox{Irr}}

\newcommand{\sltwo}{\mathit{Sl}_2}
\newcommand{\slthree}{\mathit{Sl}_3}
\newcommand{\slfour}{\mathit{Sl}_4}

\newcommand{\spfour}{\mathit{Sp}_4}
\newcommand{\Gr}{\mathcal{G}}

\newcommand{\polytinv}{\mathbb{C}[t,t^{-1}]}
\newcommand{\polyt}{\mathbb{C}[t]}
\newcommand{\IH}{\mbox{IH}}

\newcommand{\alg}{\mathcal{A}}

\newcommand{\startproof}{\noindent \textbf{Proof. }}
\newcommand{\stopproof}{$\Box$ \newline \noindent}

\newcommand{\weights}{X^*(\geoT)}

\newcommand{\rvs}{{\sf t}_{\real}^*}
\newcommand{\sg}{\Lambda^{-}}
\newcommand{\poly}{\mathcal {MV}}
\newcommand{\param}{\mathbb{B}}

\newcommand{\unip}{N}
\newcommand{\weylgon}[1]{\mbox{conv}(W\cdot{#1})}

\newcommand{\mm}{\boldsymbol \mu}

\newlength{\posttable}
\title{A Polytope Combinatorics \linebreak~for Semisimple Groups}
\author{Jared E. Anderson}
\address{Department of Mathematics, University of Massachusetts, Amherst}
\email{anderson@math.umass.edu}

\begin{document}
\begin{abstract}  

Mirkovi{\'c} and Vilonen discovered a canonical basis of algebraic cycles 
for the intersection homology of (the closures of the strata of) 
the loop Grassmannian.  The moment map images of these varieties 
are a collection of polytopes, and they may be used to compute 
weight multiplicities and tensor product multiplicities for representations 
of a semisimple group.  
The polytopes are explicitly described for a few 
low rank groups.

\vspace{3mm}

(Mathematics subject classification numbers: 14L99, 20G05)
\end{abstract}

\maketitle

%
%

\section{Introduction}

Starting with a semisimple algebraic group, we construct a collection 
of polytopes.  The central result is a method that uses these to 
decompose the tensor product of two irreducible representations.  
Each tensor product multiplicity (Littlewood-Richardson number) is 
the number of polytopes in a certain set.

The method is based on the geometry of the loop Grassmannian, 
and builds directly on the work of Mirkovi{\'c} and Vilonen~\cite{MV}.  
But all algebraic geometry is deferred until section~\ref{sec:geometry} 
since we may state the main result without it (Theorem~\ref{prop:statement}).  
We do this in section~\ref{sec:statement}, and 
follow with a lot of examples in sections~\ref{sec:examples} and~\ref{sec:examples2}.
Much may be gained from just these first sections without ever 
understanding what the loop Grassmannian is.

Some background in geometry in representation theory is discussed 
in section~\ref{sec:geometry}.  
One begins by fixing an algebraic group, 
and constructing from it a space, the loop Grassmannian.  
The representation theory of the (Langlands dual) group 
is known to be closely related to the geometry of the loop Grassmannian.  
This relationship was made more explicit with 
Mirkovi{\'c} and Vilonen's discovery of a collection of 
singular algebraic varieties in the loop Grassmannian, 
which we call MV-cycles.  
In terms of geometry, they provide a canonical basis 
for the intersection homology of the closure of 
each stratum of the loop Grassmannian.  
In terms of representation theory, they provide a canonical basis 
for each irreducible representation of the group.  

Section~\ref{sec:momentmap} provides 
a definition of the polytopes as moment map images of MV-cycles.  
The rest of the paper consists mainly of the proof of Theorem~\ref{prop:statement}.   
Section~\ref{sec:fibers} contains the main geometric idea 
 and is of interest in its own right (Theorem~\ref{prop:fiber}).  
The last section provides a glimpse at a 
closely related Hopf algebra.

\section{Statement of Results}
\label{sec:statement}

Let $\repG$ be a connected semisimple complex algebraic group of rank $n$.  
Choose a maximal torus $\repT\subset \repG$.  
We will be considering polytopes in the real 
$n$-dimensional vector space in which pictures of roots and weights 
are usually drawn:  
this is the dual $\rvs$ of the Lie algebra of the split real form of $\repT$.

We will define a collection of polytopes, $\poly=(P_\phi)_{\phi\in\param}$ in $\rvs$
(Definition~\ref{def:MVpolytopes} on page~\pageref{def:MVpolytopes}).  
Let $R^-$ denote the set of negative roots and   
$\sg$ the semigroup they generate.
The parameter set $\param$ is naturally graded by $\sg$: $\param=\bigcup_{\nu\in\sg}\param_\nu$.
It turns out that this parameter set is not important in the theorem below, 
since the parametrization $\phi\mapsto P_\phi$ is injective~\cite{AM}; 
so we will really be counting polytopes.  
We also write $\poly_\nu=\{P_\phi|\phi\in \param_\nu\}$ 
so that $\poly=\bigcup_{\nu\in\sg}\poly_\nu$.

Among our polytopes will be (shifts of) those familiar from representation 
\mbox{theory}: the convex hull of the weights in an irreducible representation of $\repG$; 
$\weylgon{\lambda}$ denotes the convex hull of the 
Weyl group orbit through a weight $\lambda$.  

\begin{theorem}
Weight multiplicities and tensor product multiplicities may be calculated according to the following 
rules.

\begin{enumerate}

\item If $V_\lambda$ is an irreducible representation of $\repG$ with dominant weight 
$\lambda$, then the multiplicity of weight $\nu$ 
in $V_\lambda$ equals the number of $\phi\in\param_{\nu-\lambda}$ for which 
$P_\phi+\lambda\subseteq \weylgon{\lambda}$.

\item If $V_\lambda$ and $V_\mu$ are irreducible representations of $\repG$ 
with dominant weights $\lambda$ and $\mu$, and $\nu$ is any dominant weight, 
then the multiplicity of $V_\nu$ in $V_\lambda \otimes V_\mu$ 
equals the number of $\phi\in \param_{\nu-\mu-\lambda}$ for which 
$P_\phi+\lambda\subseteq \weylgon{\lambda}\cap(\weylgon{{-\mu}}+\nu)$.

\end{enumerate}

\noindent

\label{prop:statement}
\end{theorem}

The polytopes in $\poly$ are called MV-polytopes.  
The ``MV'' stands for Mirkovi{\'c} and Vilonen, who discovered a 
collection of algebraic varieties called MV-cycles.  
MV-polytopes will be defined as moment map images of MV-cycles.  
Part 1 of the theorem is little more than a translation of some 
of the algebraic geometry of  Mirkovi{\'c} and Vilonen into the language of 
polytopes.  Part 2, however, depends on some more 
geometry described in section~\ref{sec:fibers}.

\section{Examples of Polytopes}
\label{sec:examples}

We explicitly describe the collection of polytopes $\poly$ for a few low rank groups: 
$\sltwo$, $\slthree$, $\spfour$, $\slfour$.  This is not known for other groups, although 
Mirkovi{\'c} and I have 
a conjecture that inductively constructs the polytopes for any group~\cite{AM}.

Before describing {\it all} the polytopes for these groups, we introduce them with a picture 
of eight of them, which count the weight multiplicities in the adjoint representation of $\slthree$.  
The weight multiplicity is 1 at each outer vertex, and 2 at the central vertex.

{\par 
\begin{center}
\includegraphics[width=12cm]{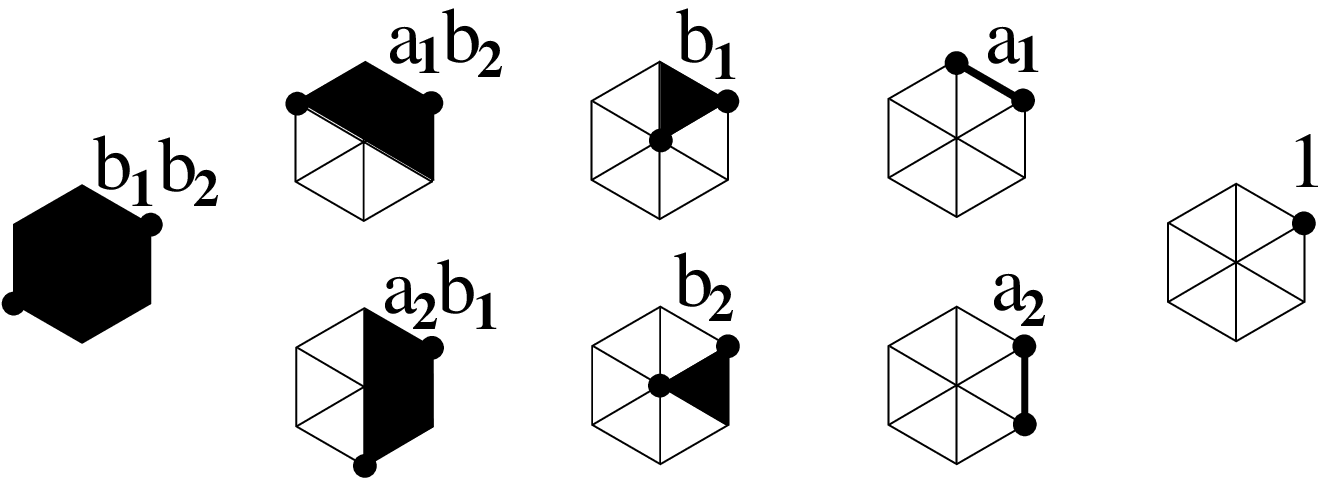}
\end{center}
}

To describe the polytopes, it is easiest to first introduce a commutative algebra $\alg$ with basis $\poly$.  
This algebra is discussed briefly in section~\ref{sec:hopf}; 
conjecturally it is the algebra 
of functions on the unipotent radical 
of a Borel subgroup~\cite{A}.  
In these examples we give ad hoc definitions of $\alg$ by 
a natural set of generators and relations.  
Then a collection of monomials in the generators is defined.  
A polytope is associated to each generator, 
and then to each monomial 
by taking the Minkowski sum of the factors.  
(The Minkowski sum of 
two sets $A$ and $B$ is the 
set of sums $\{a+b|a\in A, b\in B\}$.)

\subsection{$\sltwo$}
$\alg=\mathbb{Z}[a]$ and the MV-polytopes $[k\alpha,0]$ correspond to monomials $a^k$.   
($\alpha$ is the negative root in $\rvs=\real$.)

\subsection{$\slthree$}

The algebra $\alg$ has four generators:
{\par 
\begin{center} 
\includegraphics[width=12cm]{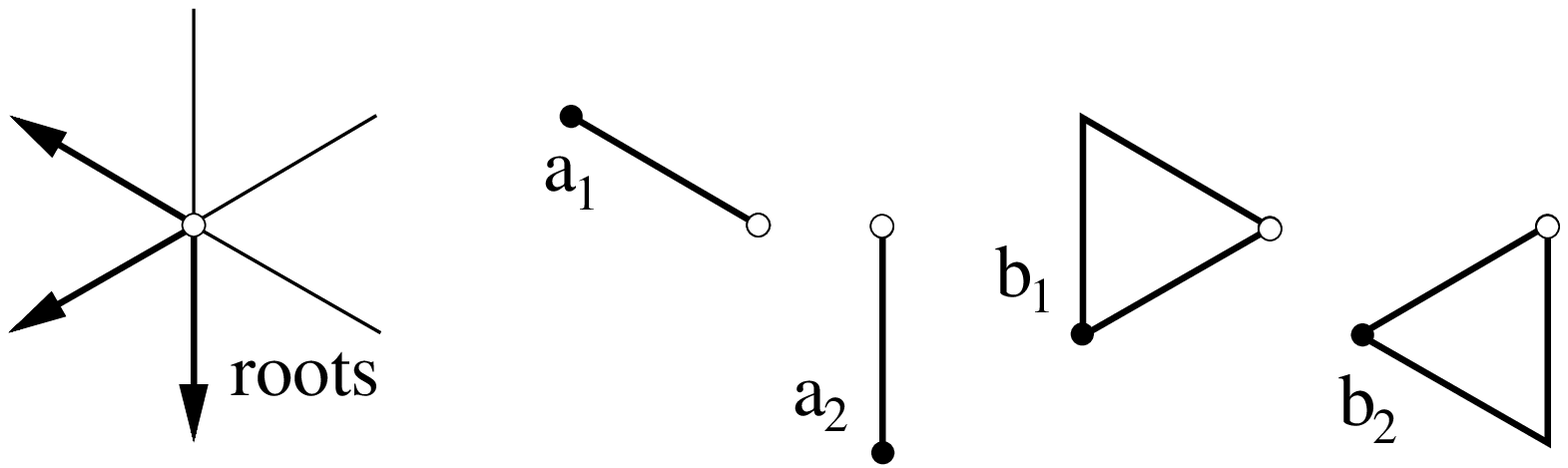}
\end{center}
}

\noindent
(The negative roots are indicated.  
For each polytope, the vertex with highest weight is at the origin.)  
There is a single relation, $a_1a_2=b_1+b_2$, which we put in a diagram:
{\par 
\begin{center} 
\includegraphics[width=10cm]{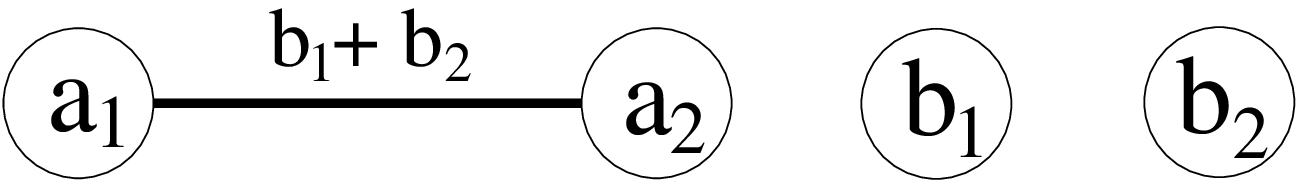}
\end{center}
}
\noindent
The MV-polytopes correspond to monomials of the form $a_1^ib_1^jb_2^k$ 
or $a_2^ib_1^jb_2^k$, i.e. those monomials in the generators such that 
not both $a_1$ and $a_2$ occur.  These 
are the monomials that cannot be further simplified using the relations.
Again, the MV-polytopes are 
found by taking Minkowski sums.  For example, 
the regular hexagon $b_1 b_2$ is the Minkowski sum of the two triangles.
In general, the monomials $b_1^j b_2^k$ give shifts of the 
symmetric hexagons $\weylgon{\lambda}$. One can think of an 
arbitrary MV-polytope as one of these hexagons stretched some length 
in the direction of either $a_1$ or $a_2$.

As an aside, note that the relation $a_1a_2=b_1+b_2$ has an interpretation in 
terms of Minkowski sum:  the Minkowski sum of the two line segments $a_1$ and $a_2$ 
is a parallelogram which equals the union of the two triangles $b_1$ and $b_2$.  
In general, if $\alpha \beta=\sum_{\gamma\in\poly} n_\gamma \gamma$ then 
the Minkowski sum of $\alpha$ and $\beta$ equals the union of the $\gamma$ for which 
$n_\gamma \neq 0$. (See~\cite{A}.)

The above definition places 
the highest weight vertex of $P\in\poly$ at the origin.
Then the lowest weight vertex $\nu$ identifies the graded 
piece $\poly_\nu$ in which $P$ lies.  

\subsection{$\spfour$}

The algebra $\alg$ has eight generators:
{\par 
\begin{center} 
\includegraphics[width=12cm]{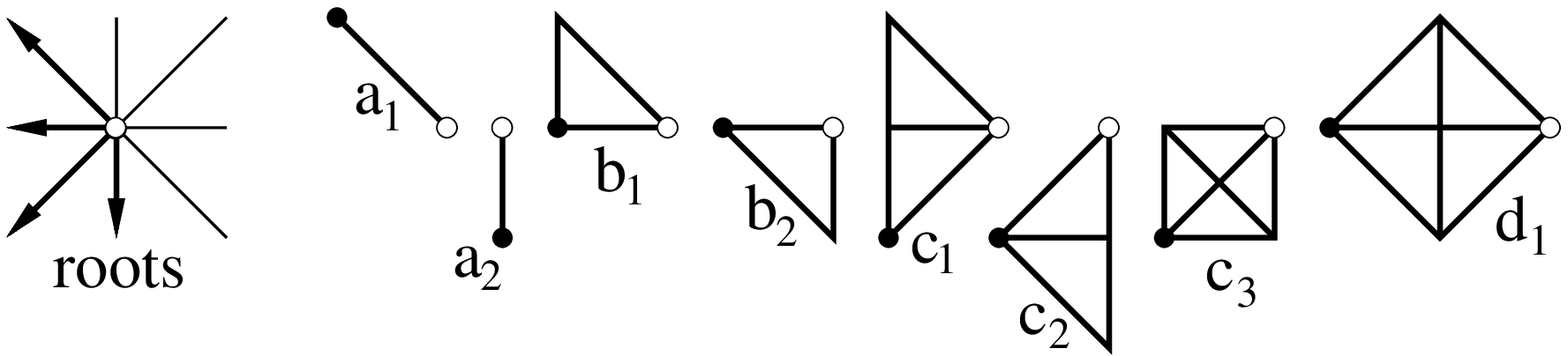}
\end{center}
}

\noindent
(Again the negative roots indicated, and the highest weight vertex 
is always at the origin.)  
There are nine relations:

{\par 
\begin{center} 
\includegraphics[width=10cm]{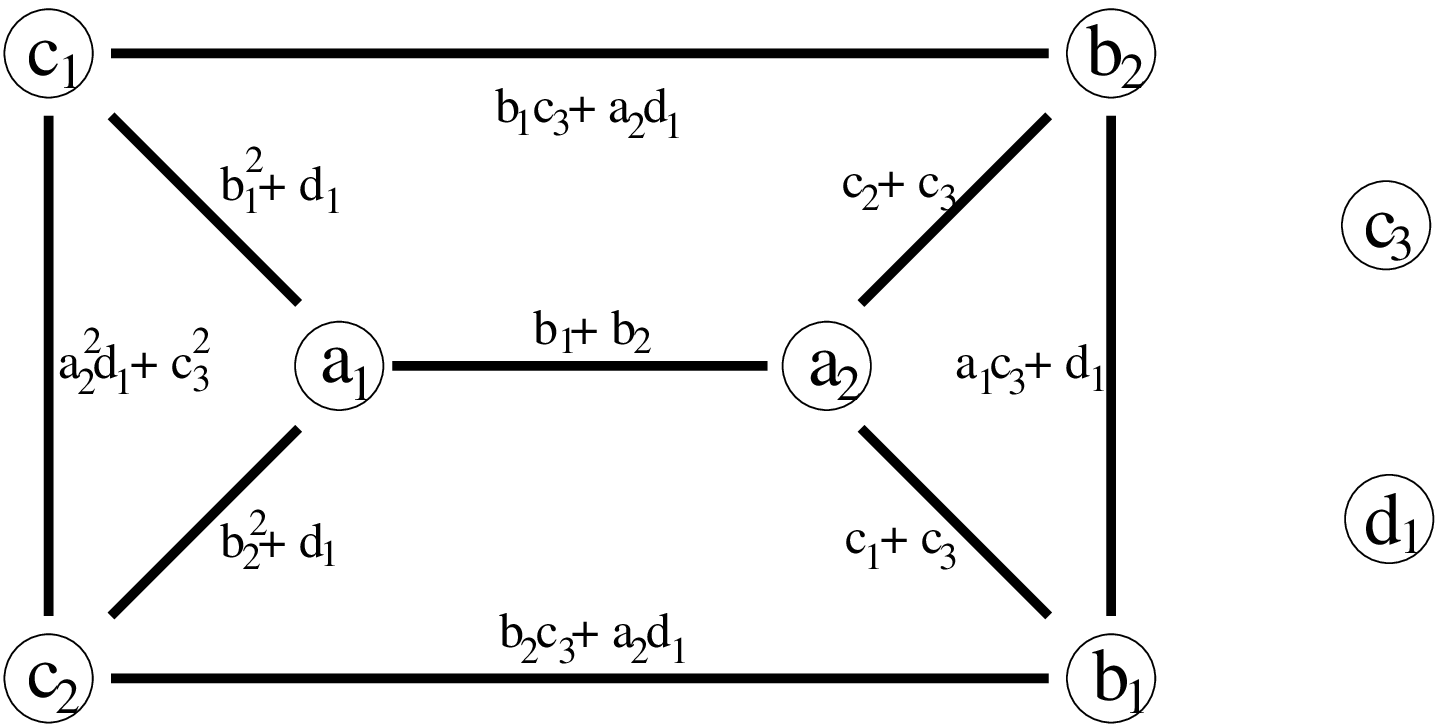}
\end{center}
}

\noindent
The MV-polytopes correspond to monomials of these forms: 
$a_1^i b_1^j c_3^k d_1^l$,
$a_1^i b_2^j c_3^k d_1^l$,
$a_2^i c_1^j c_3^k d_1^l$,
$a_2^i c_2^j c_3^k d_1^l$,
$b_1^i c_1^j c_3^k d_1^l$,
$b_2^i c_2^j c_3^k d_1^l$.  
Equivalently, these are all monomials in the generators, 
no two factors of which are joined by a line in the diagram of 
relations.  
The monomials $c_3^k d_1^l$ give shifts of the symmetric octagons 
$\weylgon{\lambda}$.  
The grading is specified exactly as 
it was for $\slthree$.

\subsection{$\slfour$}
The example of $\slfour$ is very rich, and harder to think about since 
the polytopes are three-dimensional.
There are twelve generators:  
{\par 
\begin{center} 
\includegraphics[width=12cm]{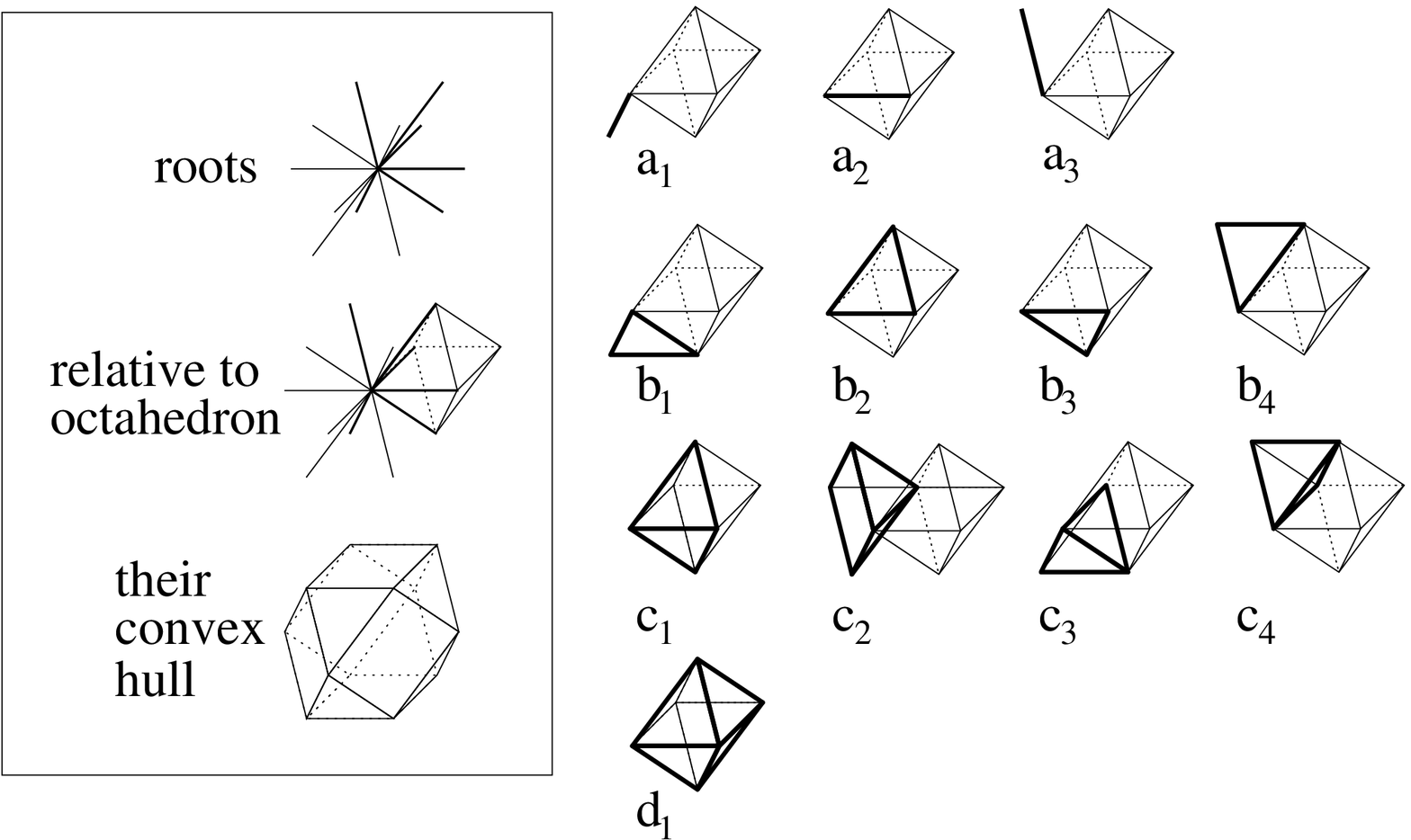}
\end{center}
}
\noindent
These are drawn relative to an octahedron    
whose front left vertex is the origin of $\rvs=\real^3$.  
The three line segments $a_1,a_2,a_3$ identify
the simple negative roots.  
The $b_i$ are triangles; $c_1,c_2$ are square-based pyramids; 
$c_3,c_4$ are tetrahedra; $d_1$ is an octahedron.  
Here are the fifteen relations: 
{\par 
\begin{center} 
\includegraphics[width=12cm]{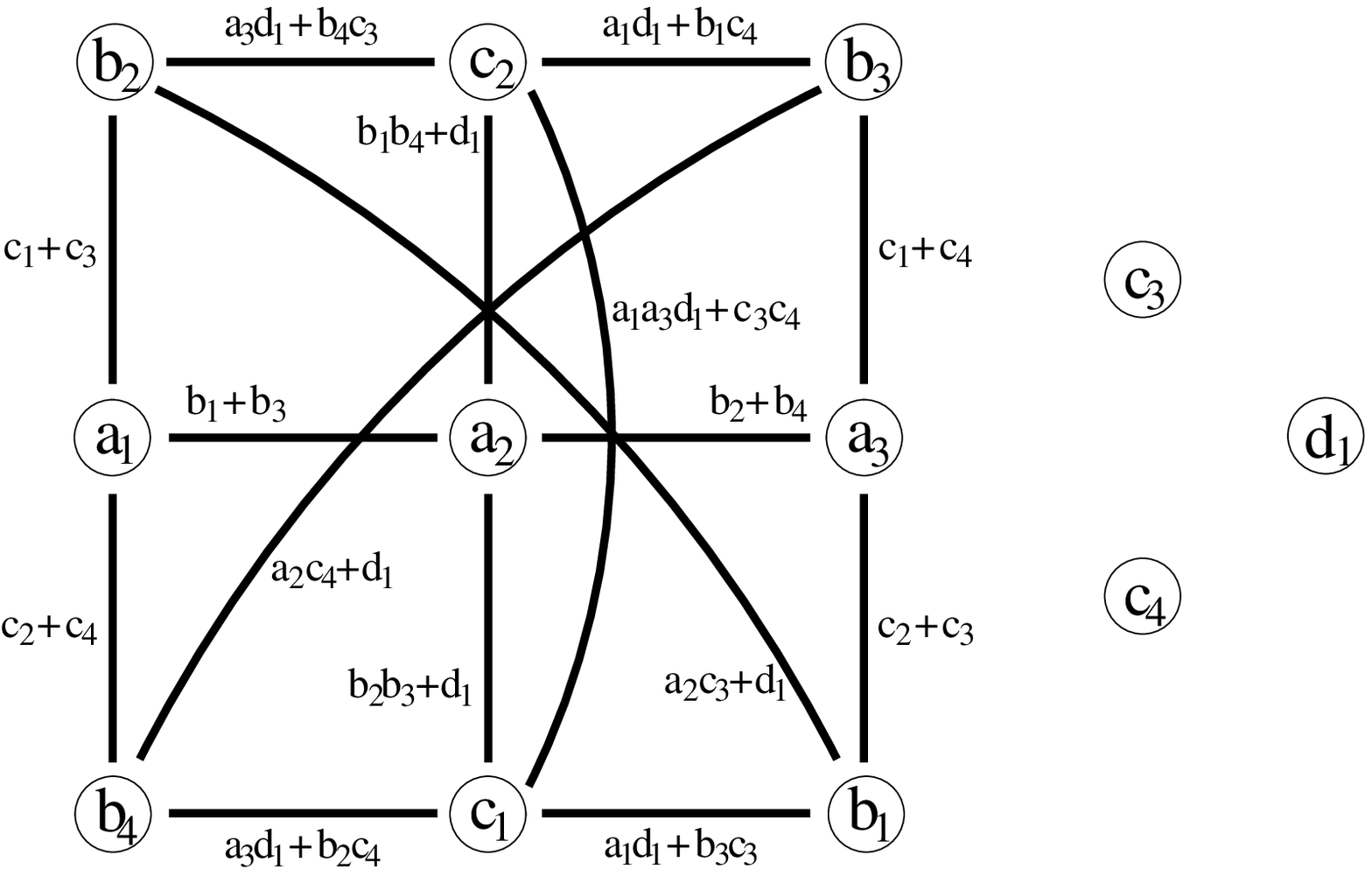}
\end{center}
}
\noindent
As before, the collection of monomials consists 
of those for which no two factors are joined by a line in this 
diagram.  

\section{Examples of Multiplicity Calculations}
\label{sec:examples2}

We use the MV-polytopes for $\spfour$ in two examples, 
illustrating the two parts of Theorem~\ref{prop:statement}.  

\subsection{Example: $\spfour$ weight multiplicity calculation}

This illustrates part 1 of Theorem~\ref{prop:statement}.  
The octagon below corresponds to an irreducible 
representation of $\spfour$.  
Suppose 
we want to know the weight multiplicity at the 
indicated weight $\nu$.  
$\poly_\nu$ consists of the six polygons pictured below.  
Each has one 
vertex, distinguished by a black dot, which must be 
placed at the weight $\nu$.   
Five of the polygons are contained in the octagon, 
and one is not (its bottom right corner sticks out).  
So the weight multiplicity is 5.
{\par 
\begin{center} 
\includegraphics[width=12cm]{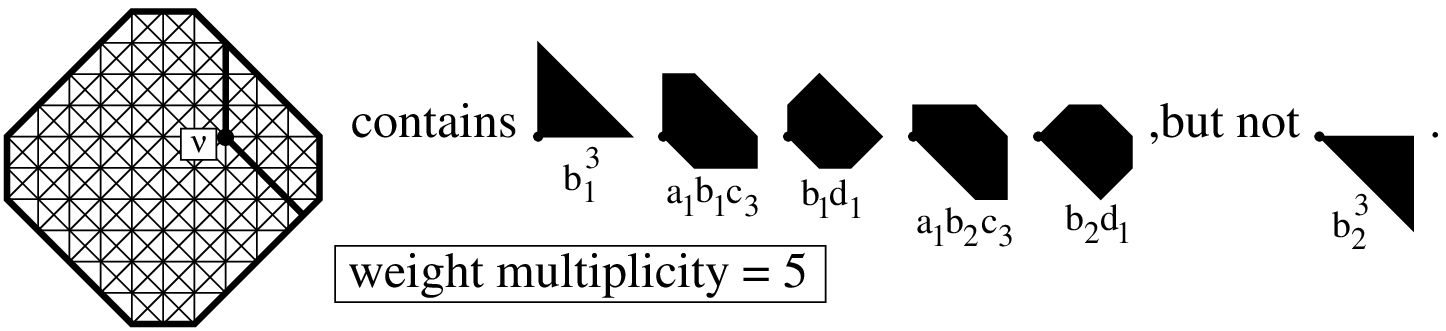}
\end{center}
}

\subsection{Example: $\spfour$ tensor product multiplicity calculation}

This illustrates part 2 of Theorem~\ref{prop:statement}.  
Suppose we want to decompose the tensor product of the two 
irreducible $\spfour$ representations indicated in the picture.
{\par 
\begin{center}
\includegraphics[width=12cm]{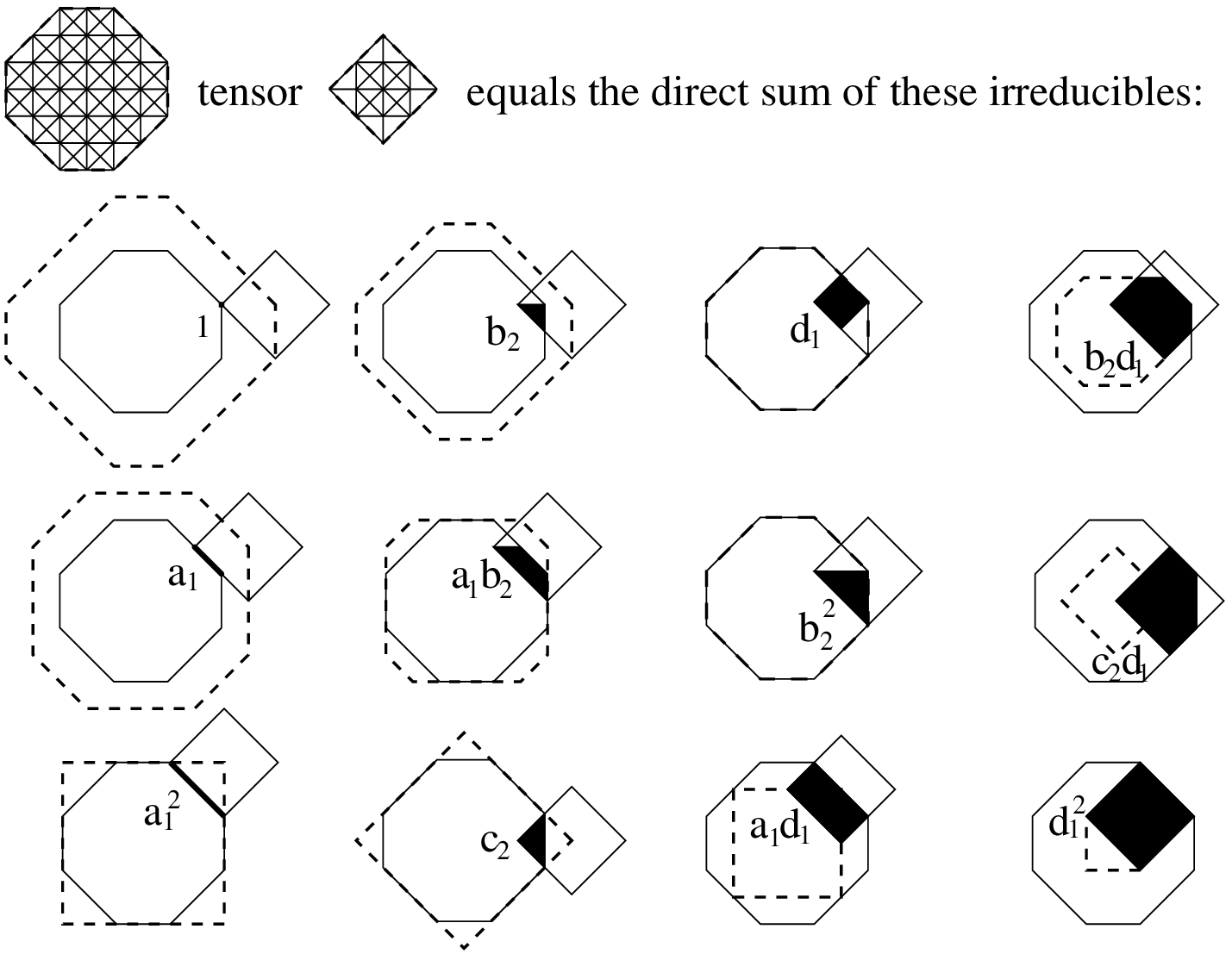} 
\end{center}
}
\noindent
The dashed lines are the outlines of the irreducible summands 
that occur with nonzero multiplicity.    
The MV-polygons that were counted to find the multiplicities 
are drawn in black; 
these 
were the ones contained in the intersection of 
the big octagon $\weylgon{\lambda}$ and the translated square 
$\weylgon{{-\mu}}+\nu$.  
Each summand has multiplicity 1, except for 
one which has multiplicity 2: $d_1, b_2^2$.

\subsection*{Remark}

This calculation has a similar flavour to the convolution 
of the characters of the two representations, where one computes 
a sum of products of weight multiplicities instead of counting 
polytopes.  
After convolving the characters, one uses the Weyl character formula 
to find the irreducible summands.  
Our method required no knowledge of weight multiplicities, 
no arithmetic, and no separate calculation to extract the 
irreducible summands.  
Of course, saying that 
this method is more efficient than using the Weyl character 
formula is not saying very much.  

One might hope that this would provide an explicit 
formula for the Littlewood-Richardson numbers, 
but it seems only to provide an algorithm.

\section{Geometry in Representation Theory}
\label{sec:geometry}

\subsection{Loop Grassmannian}

Our basic object of study is the loop Grassmannian (affine Grassmannian), which is 
an infinite dimensional space associated to a complex algebraic group $\geoG$, 
taken to be connected and semisimple.   
It is the quotient of groups $\Gr=\geoG(\mathcal{K})/\geoG(\mathcal{O})$ 
where $\mathcal{O}=\fps$, the ring of formal power
series, and $\mathcal{K}=\fls$, its field of fractions, the ring of formal Laurent series.

One can show that if we use the rings $\polytinv$ and $\polyt$ instead
of $\mathcal{K}$ and $\mathcal{O}$ in this definition, then we get the
same set.  Since the numerator $\geoG(\polytinv)$ is a group of maps from
the unit circle to $\geoG$ (let $t=e^{i\theta}$), we understand the use of
the word {\it loop}.  There is a model of $\Gr$ due to Lusztig which
describes it as a set of subspaces of an infinite dimensional vector
space, which explains the use of {\it Grassmannian}~\cite{L}.  

$\Gr$ may be
realized as an increasing union of finite dimensional complex
projective varieties by filtering $\geoG(\mathcal{K})$ by
order of pole~\cite{L}.  
The orbits of the  
action of $\geoG(\mathcal{O})$ on $\Gr$ by left multiplication 
provide a stratification of $\Gr$, and are indexed
by the dominant coweights of $\geoG$.  Here is
our notation for this:  Fix a maximal torus $\geoT\subseteq \geoG$.  
Any coweight $\lambda\in \mbox{Hom}(\cplx^*,\geoT)$ 
defines an element of
$\geoT(\mathcal{K})\subseteq \geoG(\mathcal{K})$, and hence one of $\Gr$,
which we denote $\underline{\lambda}$.  Let
$\Gr_\lambda=\geoG(\mathcal{O})\underline{\lambda}$. 
Any point in $\Gr$ is in some $\Gr_\lambda$ and
this $\lambda$ is determined up to the action of the Weyl group.

Each stratum $\Gr_\lambda$
is a complex vector bundle over a flag manifold for $\geoG$.  For a
coweight $\lambda$ in the interior of the positive Weyl chamber
this is the full flag manifold, but if $\lambda$ lies on a Weyl chamber
wall then it is some partial flag manifold.  
The closure $\overline{\Gr_\lambda}$ of a stratum consists of the union 
of all $\Gr_\mu$ with $\mu\gleq\lambda$, where $\lambda$ and $\mu$ are dominant.  
($\mu\gleq\lambda$ means that $\lambda-\mu$ is in the 
semigroup generated by the positive coroots.)
$\overline{\Gr_\lambda}$ 
is a finite dimensional complex projective algebraic variety, 
and is almost always singular.  
(See~\cite{PS}.)

\subsection{Relation to representation theory}

It turns out that what gives information about representation theory
is the intersection homology of Goresky and
MacPherson~\cite{GM1,GM2}.  This is a homology theory defined on a
stratified space, using cycles 
that satisfy bounds on the allowed dimensions of intersections with 
lower dimensional strata.  
Work of Drinfeld, Ginsburg, Lusztig, 
Mirkovi{\'c}, and Vilonen \cite{BD,VG,MV} 
shows that the intersection homology $\IH(\overline{\Gr_\lambda})$ 
of the closure of a stratum 
may be regarded as the vector 
space underlying the representation $V_\lambda$ of the Langlands 
dual group $\geoG^L=\repG$.  
(Here $\lambda$ is both a coweight of $\geoG$ and a
weight of $\repG$.)
They actually show much more: the category of 
$\geoG(\mathcal{O})$-equivariant perverse sheaves on $\Gr$ is given a tensor product, and is  
equivalent to the tensor category of representations of $\repG$.

\subsection{MV-cycles}
\label{sec:MV}

Mirkovi{\'c} and Vilonen discovered a canonical basis of algebraic
cycles for the intersection homology of stratum closures
$\overline{\Gr_\lambda}$~\cite{MV}.  
These MV-cycles are projective varieties in the loop Grassmannian.  
The intersection homology of 
$\overline{\Gr_\lambda}$ is represented by those MV-cycles
contained in it, but not contained in any $\overline{\Gr_\mu}$ 
with $\mu\glt\lambda$.  

To define MV-cycles, we must 
first make some choices.  Fix opposite Borel subgroups of $\geoG$
that intersect in the maximal torus $\geoT$.  Denote by $\unip$ and $\unip^-$ their 
unipotent radicals.  
{\em We choose the positive roots to be roots in $\unip^-$}.
Then
for any coweights $\lambda,\mu$, we let
$S_\lambda=\unip(\mathcal{K})\underline{\lambda}$ and
$T_\mu=\unip^-(\mathcal{K})\underline{\mu}$.  Any
$\unip(\mathcal{K})$-orbit contains a unique $\underline{\lambda}$,  
as does any $\unip^-(\mathcal{K})$-orbit.
These orbits have both infinite
dimension and infinite codimension in $\Gr$.  
There are simple closure relations:  
$\overline{S_\lambda}=\bigcup_{ \xi\ggeq\lambda } S_\xi$ and 
$\overline{T_\mu}=\bigcup_{\eta\gleq\mu} T_\eta$~\cite{MV}.  

\begin{definition}
Let $\overline{\Gr_\lambda}$ be the closure of a stratum of the loop Grassmannian, 
where $\lambda$ is chosen dominant.  
Let $\nu$ be a coweight of $\geoG$ with $\underline{\nu}\in\overline{\Gr_\lambda}$.  
The MV-cycles for $\overline{\Gr_\lambda}$ at $\nu$, relative to $\unip$, 
are the irreducible components of $\overline{S_\nu \cap \Gr_\lambda}$.
Equivalently, they are those irreducible components of $\overline{S_\nu\cap T_\lambda}$ 
contained in $\overline{\Gr_\lambda}$.
\label{def:MVcycles}
\end{definition}

\noindent
The equivalence of the two definitions is a consequence of the 
dimension calculations of~\cite{MV}:  
both $\overline{S_\nu\cap T_\lambda}$ and $\overline{S_\nu \cap \Gr_\lambda}$ 
have pure (complex) dimension equal to the height of $\lambda-\nu$.

\begin{proposition}
The two definitions are equivalent.
\label{prop:equivdef}
\end{proposition}
\startproof
Suppose $a\in\Irr(\overline{S_\nu\cap T_\lambda})$ 
(i.e. $a$ is an irreducible component of this variety) 
and $a\subseteq\overline{\Gr_\lambda}$.  
Evidently $a\subseteq\overline{S_\nu}\cap\overline{\Gr_\lambda}$.  
Since $\overline{S_\nu}=\bigcup_{\xi\ggeq\nu} S_\xi$ 
and $\overline{\Gr_\lambda}=\bigcup_{\mbox{\scriptsize dom.} \eta\gleq\lambda} \Gr_\eta$, 
we see that $\overline{S_\nu}\cap\overline{\Gr_\lambda}=\overline{S_\nu\cap\Gr_\lambda}\cup X$ 
where every component of $X$ has dimension strictly less than 
height$(\lambda-\nu)$, the dimension of $a$.  
So $a\in\Irr(\overline{S_\nu\cap \Gr_\lambda})$.  

Conversely, suppose  $a\in\Irr(\overline{S_\nu\cap \Gr_\lambda})$.  
According to~\cite{MV}, $\ddim (\Gr_\lambda\cap T_\lambda)=\ddim(\Gr_\lambda)$; 
therefore $\overline{\Gr_\lambda}\subseteq \overline{T_\lambda}$.  
So $a\subseteq\overline{S_\nu}\cap\overline{T_\lambda}$.  
Since $\overline{S_\nu}=\bigcup_{\xi\ggeq\nu} S_\xi$ 
and $\overline{T_\lambda}=\bigcup_{\eta\gleq\lambda} T_\eta$, 
we see that $\overline{S_\nu}\cap\overline{T_\lambda}=\overline{S_\nu\cap T_\lambda}\cup Y$ 
where every component of $Y$ has dimension strictly less than 
height$(\lambda-\nu)$, the dimension of $a$.  
(In fact, I believe $Y$ is empty.)  
So $a\in\Irr(\overline{S_\nu\cap T_\lambda})$ 
and clearly $a\subseteq\overline{\Gr_\lambda}$.  
\stopproof

When $\nu=\lambda$ 
there is one MV-cycle, the point $\underline{\lambda}$.   
When $\nu=w_0 \lambda$ (where $w_0$ is the longest element of the Weyl group), 
there is also one MV-cycle, 
the whole variety $\overline{\Gr_\lambda}$.  
In general, the number of irreducible components of $\overline{S_\nu \cap \Gr_\lambda}$ 
is the dimension of the 
weight space $\nu$ in the irreducible representation $V_\lambda$~\cite{MV}.

\section{The Moment Map and MV-polytopes}

\label{sec:momentmap}

The moment map $\mm$,  
for the action of the torus $\geoT$ on $\Gr$, 
is a map from 
$\Gr$ to
$\lie(\geoT)^*$. 
We use the Killing form to identify $\lie(\geoT)^*$ with $\lie(\geoT)$, 
which in turn is canonically identified with $\lie(\repT)^*$.  
We will view $\lie(\repT)^*$ as the codomain of $\mm$,   
since moment map images will have to do 
with the representation theory of $\repG$, which 
is usually pictured in the real subspace $\rvs$ of $\lie(\repT)^*$. 

We define 
$\mm$ as the restriction of the moment map on a projective space.  
Let $\mathcal L$ be an ample line bundle on $\Gr$ and 
$\Gamma(\Gr, \mathcal{L})$ the vector space of global sections.  
Then $\Gr$ naturally embeds in the projective space 
$\mathbb{P}(V)$ where $V=\Gamma(\Gr,\mathcal{L})^*$ 
by 
mapping 
$x\in\Gr$ to the point 
determined by the line in V dual to the 
hyperplane $\{s\in\Gamma(\Gr,\mathcal{L}) | s(x)=0 \}$.
The action of the torus $\geoT$ on V decomposes it into 
eigenspaces: $V=\bigoplus_{\nu\in\weights} V_\nu$, where $\weights$ denotes the weights.
Choose an inner product on $V$ that is invariant under the action of 
the maximal compact subgroup of $T$, so that this decomposition is orthogonal.  
Then, given $v = \sum 
v_\nu \in V$, we define 
$\mu([v])=\sum \frac{|v_\nu|^2}{|v|^2} \nu$, 
the usual moment map on a projective space.

\begin{proposition} \

\begin{enumerate}

\item The fixed points of the torus action are the  
$\underline{\nu}\in\Gr$  ($\nu$ a coweight of $\geoG$); 
and $\mm(\underline{\nu})=\nu$.

\item If $X$ is a one-dimensional torus orbit then $\mm(X)$ 
is a line segment in a root direction joining two weights.  

\item If $a$ is an MV-cycle then $\mm(a)$ is a convex polytope.  
Its vertices are a collection of weights parametrized by the Weyl group, 
possibly in a degenerate way.

\item $\mm(\overline{\Gr_\lambda})=\weylgon{\lambda}$

\item If $X\subseteq\Gr$ is any compact torus-invariant variety and $\eta$ is any 
coweight then $\mm(\eta X)$=$\eta+\mm(X)$.
\end{enumerate}

\end{proposition}
\startproof
(1)~It is easy to check that any $\underline{\nu}$ is a fixed point.  
That there aren't others follows, for instance, from decomposing $\Gr$ into 
$\unip(\mathcal{K})$ orbits.  For a coweight $\nu$, we have $v=v_\nu$ so that $\mm(\underline{\nu})=\nu$.
(2)~$\mm(X)$ is certainly a line segment joining two weights~\cite{GM3}.  
That it lies in a root direction follows from viewing $\Gr$ as a partial 
flag variety and knowing the $\geoT$-invariant curves in a flag variety~\cite{C}.  
(3)~The moment map image of any compact irreducible torus-invariant 
variety is the convex hull of the images of its $\geoT$-fixed 
points~\cite{B,GM3}.  How the vertices are parametrized 
by the Weyl group is explained in~\cite{A}; this fact 
is not used here.  
(4)~$\mm(\overline{\Gr_\lambda})$ is the convex hull of the images 
of its fixed points.  
Since the fixed points of $\Gr_\lambda$ are the $W\cdot\underline{\lambda}$, 
we have $\mm(\overline{\Gr_\lambda})\supseteq\weylgon{\lambda}$.  
But by the closure relations for strata, all the other fixed points 
of $\overline{\Gr_\lambda}$ are also in this convex set.
(5)~For each fixed point $\underline{\xi_i}$ of $X$, we have 
$\mm(\eta \underline{\xi_i})=\mm(\underline{\eta+\xi_i})=\eta+\xi_i=\eta+\mm(\underline{\xi_i})$.  
The $\eta\underline{\xi_i}$ are the fixed points of $\eta X$.  
The statement follows since $\mm(X)$ and $\mm(\eta X)$ are 
the convex hulls of the $\mm(\underline{\xi_i})$ and $\mm(\eta \underline{\xi_i})$ respectively.
\stopproof

\begin{definition}
For each $\nu\in\sg$, let
$\param_\nu=\Irr(\overline{S_\nu\cap T_0})$,   
and for each irreducible component $\phi$ in $\param_\nu$, let 
$P_\phi$ be its moment map image $\mm(\phi)$.  
So $\poly=\bigcup_\nu \poly_\nu$ 
where $\poly_\nu=\{\mm(a)|a\in \param_\nu \}$.

\label{def:MVpolytopes}
\end{definition}

\section{Proof of Weight Multiplicity Calculation \\(part 1 of Theorem 1)}

\begin{proposition}
Let $\underline{\xi}$ be any fixed point in the closure 
of the stratum $\overline{\Gr_\eta}$, where $\eta$ is chosen dominant.  
An irreducible component $a$ of $\overline{S_\xi\cap T_\eta}$   
is an MV-cycle if and only if its moment map image $\mm(a)$
is contained in $\mm(\overline{\Gr_\eta})$.
\label{prop:weights}
\end{proposition}
\startproof
If $a$ is an MV-cycle then it is a component of 
$\overline{S_\xi\cap\Gr_\eta}$, which is contained in $\overline{\Gr_\eta}$.  
So $\mm(a)$ is contained in $\mm(\overline{\Gr_\eta})$.

To see the other direction, we need to use 
a larger torus action of $\cplx^*\times \geoT$ on 
$\Gr=\geoG(\mathcal{K})/\geoG(\mathcal{O})$: 
the first factor acts on the
indeterminate that occurs in $\mathcal{K}$, and the $\geoT$ acts just as before.
All of the varieties in question are preserved 
by this action.  
Let us assume that $\mm(a)\subseteq\mm(\overline{\Gr_\eta})$.  
Suppose $x\in a$.  We know that $x$ 
is contained in some stratum, say $\Gr_\epsilon$.  Like all strata, 
this is a vector bundle over a flag manifold, $\geoG(\cplx)\underline{\epsilon}$, 
sometimes called the core.  
The action of
small $t\in\cplx^*$ (the first factor in $\cplx^*\times \geoT$)
on $x$ sends it close to the core.  Choosing a 
sequence of such $t_n$ converging to $0$, we can construct 
a point $y=\lim t_nx$ contained in both $a$ and the core.  
Then acting by the second factor $\geoT$ will allow us to move $y$ 
arbitrarily close to some fixed point in the core; such points 
are the Weyl translates of $\underline{\epsilon}$.  Since $a$ is closed 
and preserved 
by the torus, this fixed point is contained in $a$. 
By assumption, it follows that $\epsilon$ lies in 
$\mm(\overline{\Gr_\eta})$.  
By the closure relations for strata, this implies 
$\Gr_\epsilon\subseteq\overline{\Gr_\eta}$.  
Hence $x\in\overline{\Gr_\eta}$ and we have $a\subseteq\overline{\Gr_\eta}$.  
So $a$ is an MV-cycle for $\Gr_\eta$.
\stopproof

\noindent
{\bf Proof of Part 1 of Theorem~1. }
According to~\cite{MV}, the weight multiplicity at weight $\nu$ 
in an irreducible representation $V_\lambda$ 
equals the number of MV-cycles at weight $\nu$, i.e. the number of 
components of $\overline{S_\nu\cap T_\lambda}$ that are contained in 
$\overline{\Gr_\lambda}$.    
By the preceding proposition, this is the number of 
$a\in\Irr(\overline{S_\nu\cap T_\lambda})$ such that 
$\mm(a)$ is contained in $\mm(\overline{\Gr_\lambda})=\weylgon{\lambda}$.  
But this is the same as the number of 
$a\in\Irr(\overline{S_{\nu-\lambda}\cap T_0})$  
such that $\mm(a)+\lambda$ is contained in $\weylgon{\lambda}$.
Therefore the weight multiplicity is the number of 
$\phi\in\param_{\nu-\lambda}$
such that $P_\phi+\lambda\subseteq \weylgon{\lambda}$. \stopproof

\section{Fibers of the Convolution Map}

\label{sec:fibers}

This section describes the geometric idea that, when translated 
into the language of polytopes, yields the tensor product 
calculation of Theorem~\ref{prop:statement}. 
Given two stratum closures $\overline{\Gr_\lambda}$, $\overline{\Gr_\mu}$ of the loop Grassmannian, 
we recall the 
construction of their twisted product ${\overline{\Gr_\lambda}} \tilde{\times} \overline{\Gr_\mu}$
and a convolution map $\pi$ 
from this to $\overline{\Gr_{\lambda+\mu}}$.  
We show that the 
relevant irreducible components of the fibers of 
this map are MV-cycles.  

\subsection{Relative position convolution}
\label{sec:relposconv}

We describe the relative position convolution of 
the closures of two 
strata in the loop Grassmannian.  Given two points $a\geoG(\mathcal{O})$
and $b\geoG(\mathcal{O})$ in $\Gr$, one can ask: what is the relative
position of $b\geoG(\mathcal{O})$ with respect to $a\geoG(\mathcal{O})$?  By
definition, this means: in what stratum is $a^{-1}b\geoG(\mathcal{O})$?
So, given closures of strata $\overline{\Gr_\lambda}$ and $\overline{\Gr_\mu}$, 
we can define an algebraic variety 
often called the twisted product:
\[
\overline{\Gr_\lambda} \tilde{\times} \overline{\Gr_\mu} =
 \{ (a\geoG(\mathcal{O}),b\geoG(\mathcal{O}))\in\overline{\Gr_\lambda}\times\overline{\Gr_{\lambda+\mu}} 
 \, | \, a^{-1}b\geoG(\mathcal{O})\in \overline{\Gr_\mu}\} 
\]
This space maps to $\overline{\Gr_{\lambda+\mu}}$ by projection $\pi$ onto the second 
factor, a stratified map of algebraic varieties.  
The decomposition theorem~\cite{M,BBD} applied to this map 
decomposes the intersection homology of 
$\overline{\Gr_\lambda} \tilde{\times} \overline{\Gr_\mu}$.  
The map is known to be semi-small~\cite{MV}, which means that the dimension 
of a fiber over a stratum $\Gr_\nu$ is not larger than 
half the codimension of the stratum in $\overline{\Gr_{\lambda+\mu}}$.  
Because of this, the decomposition theorem has a particularly simple 
form in this case:
\[
\IH(\overline{\Gr_\lambda})\otimes \IH(\overline{\Gr_\mu})
\cong
\IH({\overline{\Gr_\lambda} \tilde{\times} \overline{\Gr_\mu}}) \cong 
\bigoplus_{\Gr_\nu\subseteq \overline{\Gr_{\lambda+\mu}}} F_\nu \otimes \IH(\overline{\Gr_\nu})
\]
Here we take $F_\nu$ to be the vector space spanned by the 
fundamental classes 
of each component of the fiber over 
$\underline{\nu}\in\Gr_\nu$ that has maximum 
possible complex dimension---in this case  
the height of $\lambda+\mu-\nu$.  
These are called the relevant components since they are 
the ones that appear in the above decomposition.  
$\IH$ means the global intersection homology.  
It always has coefficients in the trivial local system since 
each stratum is simply connected~\cite{BD}.

The isomorphism on the left, due to~\cite{BD,VG,MV}, 
allows one to relate this geometry to a tensor product of representations:  
$V_\lambda \otimes V_\mu = \bigoplus
{V_\nu}^{\oplus \mbox{\scriptsize dim} F_\nu}$.

\subsection{Fibers are MV-cycles}

Since the twisted product 
$\overline{\Gr_\lambda} \tilde{\times} \overline{\Gr_\mu}$ 
is a subset of $\overline{\Gr_\lambda}\times\overline{\Gr_{\lambda+\mu}}$ 
that maps to $\overline{\Gr_{\lambda+\mu}}$ by projection onto the second 
factor, we may view any fiber as a subset of the first factor.  
Viewed as such, we show that 
the relevant components of the fiber at 
$\underline{\nu}\in\overline{\Gr_{\lambda+\mu}}$ are 
MV-cycles for $\overline{\Gr_\lambda}$.

\label{`fib'}
\begin{theorem}
For the map 
$\pi:\overline{\Gr_\lambda} \tilde{\times} \overline{\Gr_\mu} \rightarrow \overline{\Gr_{\lambda+\mu}}$, 
the relevant irreducible components of the fiber $\pi^{-1}(\underline{\nu})$, 
for $\nu\gleq\lambda+\mu$ dominant, are MV-cycles for $\overline{\Gr_\lambda}$.  
They are precisely those MV-cycles that are also contained in 
$\nu\overline{\Gr_{-\mu}}$. 
If $a$ is such a cycle 
then $\nu^{-1}a$ is an MV-cycle for $\overline{\Gr_{-\mu}}$
relative to the opposite choice of unipotent subgroup.
\label{prop:fiber}
\end{theorem}
\startproof
As mentioned above, we view the fiber as a subset of $\overline{\Gr_\lambda}$:  \begin{eqnarray*}
\pi^{-1}(\underline{\nu})=\pi^{-1}(\nu \geoG(\mathcal{O})) 
&\cong& 
\{a\geoG(\mathcal{O})\in\overline{\Gr_\lambda}|a^{-1}\nu \geoG(\mathcal{O})\in\overline{\Gr_\mu} \} \\
&=& 
\{ a\geoG(\mathcal{O})\in\overline{\Gr_\lambda}| \nu^{-1}a\geoG(\mathcal{O}) \in \overline{\Gr_{-\mu}} \}\\
&=&
\overline{\Gr_\lambda}\cap\nu\overline{\Gr_{-\mu}}
\end{eqnarray*}

\noindent
We claim that 
\begin{equation}
\overline{\Gr_\lambda}\cap \nu\overline{\Gr_{-\mu}} = (\overline{\Gr_\lambda} \cap \overline{S_{\nu-\mu}})
\cap
\nu (\overline{\Gr_{-\mu}}\cap \overline{T_{-\nu+\lambda}}).
\label{eq:MVintersection}
\end{equation}

\noindent
Obviously the left side contains the right side.  
The reverse containment is because 
$ \nu\overline{\Gr_{-\mu}}\subseteq \overline{S_{\nu-\mu}}$
and 
$ \overline{\Gr_\lambda}\subseteq \nu\overline{T_{-\nu+\lambda}}$.
Let us check the first of these; the second is similar.  
We need the closure relations for strata and for 
unipotent orbits:  
\begin{eqnarray*}
\Gr_\xi\subseteq\overline{\Gr_\eta}  \hspace{5mm}(\xi,\eta \mbox{ dominant})\hspace{8.5mm} 
& \mbox{iff} & 
\xi\gleq\eta  \hspace{4mm}\mbox{or, equivalently,} \\
\Gr_\xi\subseteq\overline{\Gr_\eta}  \hspace{5mm}(\xi,\eta \mbox{ anti-dominant}) 
& \mbox{iff} & 
\xi\ggeq\eta. \\
S_\xi\subseteq\overline{S_\eta}  \hspace{5.0mm}(\xi,\eta \mbox{ any weights})\hspace{4mm} 
& \mbox{iff} & 
\xi\ggeq\eta.
\end{eqnarray*}

\noindent
So $\overline{S_{\nu-\mu}}
= \bigcup S_{\nu-\delta}
= \bigcup\nu S_{-\delta}
= \nu \bigcup S_{-\delta} 
=\nu \overline{S_{-\mu}}$, 
each union taken over $\delta\gleq\mu$.  
Therefore it suffices to show $\overline{\Gr_{-\mu}}\subseteq\overline{S_{-\mu}}$.  
Suppose $x\in\overline{\Gr_{-\mu}}$.  
Then $x\in\Gr_{-\epsilon}$ for some $\epsilon\gleq\mu$; 
here $-\mu, -\epsilon$ are anti-dominant.  
Now, $x\in S_\delta$ for some $\underline{\delta}\in\overline{\Gr_{-\epsilon}}$
since $x$ lies in some MV-cycle for $\Gr_{-\epsilon}$ \cite{MV}; 
so $\delta\ggeq-\epsilon\ggeq-\mu$.  
Hence $S_\delta\subseteq\overline{S_{-\mu}}$,  
and $x\in\overline{S_{-\mu}}$ as required.

Now, MV-cycles for $\overline{\Gr_\lambda}$ at weight $\nu-\mu$ 
are the irreducible components of 
$\overline{\Gr_\lambda\cap S_{\nu-\mu}}$; 
similarly MV-cycles for $\overline{\Gr_{-\mu}}$ at weight $-\nu+\lambda$ 
relative to the opposite Borel are components of 
$\overline{\Gr_{-\mu}\cap T_{-\nu+\lambda}}$.  
These sets may be smaller than the sets
$\overline{\Gr_\lambda} \cap \overline{S_{\nu-\mu}}$
and 
$\overline{\Gr_{-\mu}}\cap \overline{T_{-\nu+\lambda}}$ 
in equation~\ref{eq:MVintersection}, 
but only by lower dimensional components.  
For instance, if we write 
$\overline{\Gr_\lambda}=\bigcup\Gr_{\lambda^\prime}$ 
($\lambda^\prime\gleq\lambda$ dominant) 
and 
$\overline{S_{\nu-\mu}}=\bigcup S_{\nu-\mu^\prime}$ 
($\mu^\prime\gleq\mu$), 
we see that \vspace{-2mm}
\begin{eqnarray*}
\overline{\Gr_\lambda} \cap \overline{S_{\nu-\mu}}
=
\bigcup_{\lambda^\prime,\mu^\prime} \Gr_{\lambda^\prime} \cap S_{\nu-\mu^\prime} 
= \vspace{-5cm}
\bigcup_{\lambda^\prime,\mu^\prime} \overline{\Gr_{\lambda^\prime} \cap S_{\nu-\mu^\prime}}.
\end{eqnarray*}
(To see the second equality, note that each closure in the third expression 
is a subset of the first expression.)  
The dimension of 
$\overline{\Gr_{\lambda^\prime} \cap S_{\nu-\mu^\prime}}$
is the height of $\lambda^\prime+\mu^\prime-\nu$, 
which is strictly smaller than the height of $\lambda+\mu-\nu$ 
unless $\lambda=\lambda^\prime$, $\mu=\mu^\prime$.  
Therefore $\overline{\Gr_\lambda} \cap \overline{S_{\nu-\mu}}$ 
equals $\overline{\Gr_\lambda\cap S_{\nu-\mu}}$ possibly 
together with some lower dimensional MV-cycles.   

We have shown that $\pi^{-1}(\underline{\nu})$ equals 
the intersection of 
$\overline{\Gr_\lambda \cap S_{\nu-\mu}}$
and
$\nu (\overline{\Gr_{-\mu}\cap T_{-\nu+\lambda}})$ 
possibly together with some components of dimension 
less than $\lambda+\mu-\nu$.  All statements in the 
theorem are proved.
\stopproof

\enlargethispage{1000pt}

\section{Proof of Tensor Product Multiplicity Calculation \\(part 2 of Theorem 1)}

\label{sec:tensorproducts}

The preceding section gave a geometric interpretation 
to 
the decomposition of the tensor product of two irreducible 
representations $V_\lambda$ and $V_\mu$ into irreducibles: 
the tensor product multiplicities are the dimensions of 
the $F_\nu$.  
Theorem~\ref{prop:fiber} 
suggests a method for their calculation:  
count the number of MV-cycles in each fiber.  
Part~2 of Theorem~\ref{prop:statement} describes this count 
in terms of moment map images.

For another point of view on using the loop Grassmannian 
to decompose tensor products, see~\cite{BG}, 
where a construction of Kashiwara's crystal bases is given.

\begin{proposition}
Let $\lambda$, $\mu$, $\nu$ be dominant weights 
with $\nu\gleq\lambda+\mu$.  
An irreducible component $a$ 
of $\overline{S_{\nu-\mu}\cap T_\lambda}$ 
is contained in $\overline{\Gr_\lambda}\cap\nu\overline{\Gr_{-\mu}}$ if 
and only if its moment map image $\mm(a)$ is contained in 
$\mm(\overline{\Gr_\lambda})\cap(\mm(\overline{\Gr_{-\mu}})+\nu)$.
\label{prop:tensorcount}
\end{proposition}
\startproof
Obviously containment of the varieties implies containment of the
moment map \mbox{images}.  
Conversely, 
if $\mm(a)\subseteq\mm(\overline{\Gr_\lambda})\cap(\mm(\overline{\Gr_{-\mu}})+\nu)$, 
then 
$\mm(a)\subseteq \mm(\overline{\Gr_\lambda})$ 
with 
$a\subseteq\overline{S_{\nu-\mu}\cap T_\lambda}$, 
and 
$\mm(\nu^{-1}a)\subseteq \mm(\overline{\Gr_{-\mu}})$ 
with 
$\nu^{-1}a\subseteq \overline{S_{-\mu}\cap T_{-\nu+\lambda}}$.  
Proposition~\ref{prop:weights} applied twice (once for the opposite 
unipotent) implies that $a$ is 
contained in 
$\overline{\Gr_\lambda}$ and in $\nu\overline{\Gr_{-\mu}}$ as required.
\stopproof

\noindent
{\bf Proof of Part 2 of Theorem~1. }
By Proposition~\ref{prop:fiber}, 
the multiplicity with which $V_\nu$ occurs 
in $V_\lambda\otimes V_\mu$ equals the number of MV-cycles 
for $\overline{\Gr_\lambda}$ contained in $\nu\overline{\Gr_{-\mu}}$.  
This is the number of irreducible components of 
$\overline{S_{\nu-\mu}\cap T_\lambda}$ 
contained in 
$\overline{\Gr_\lambda}\cap\nu\overline{\Gr_{-\mu}}$.  
By Proposition~\ref{prop:tensorcount}, this is the number of 
$a\in\Irr(\overline{S_{\nu-\mu}\cap T_\lambda})$ such that 
$\mm(a)$ is contained in 
$\mm(\overline{\Gr_\lambda})\cap(\mm(\overline{\Gr_{-\mu}})+\nu)=\weylgon{\lambda}\cap(\weylgon{{-\mu}}+\nu)$.  
This is the same as the number of 
$a\in\Irr(\overline{S_{\nu-\mu-\lambda}\cap T_0})$ such that 
$\mm(a)+\lambda$ is contained in 
$\weylgon{\lambda}\cap(\weylgon{{-\mu}}+\nu)$.
Therefore the tensor product multiplicity is the number of 
$\phi\in\param_{\nu-\mu-\lambda}$ 
such that $P_\phi+\lambda\subseteq \weylgon{\lambda}\cap(\weylgon{{-\mu}}+\nu)$. \stopproof

\newpage
\subsection*{Remark.}

The two parts of Theorem~\ref{prop:statement} suggest that 
weight multiplicities and tensor product multiplicities are closely related.  
This is known, but since MV-cycles provide 
such a simple explanation of this,
we highlight it as a theorem:

\begin{theorem}
Suppose $\lambda$, $\mu$ are dominant weights and $\delta\ggeq 0$ is
such that $\lambda+\mu-\delta$ is dominant.
Then the multiplicity of $V_{\lambda+\mu-\delta}$ in $V_\lambda\otimes V_\mu$ 
is less than or equal to the multiplicity of the weight $\lambda-\delta$ in $V_{\lambda}$. 
By Kostant's formula, this in turn is less than or equal to the 
number of ways of writing $\delta$ as a sum of positive roots; 
moreover, this bound is 
sharp in the sense that, given $\delta$, if $\lambda$ and $\mu$ are chosen 
sufficiently large, then for all $0\gleq \epsilon\gleq\delta$ with $\lambda+\mu-\epsilon$ dominant, 
the multiplicity of $V_{\lambda+\mu-\epsilon}$ in $V_\lambda\otimes V_\mu$ exactly equals 
the number of ways of writing $\epsilon$ as a sum of positive roots.
\end{theorem}
\startproof
The multiplicity of the weight $\lambda-\delta$ in $V_{\lambda}$ equals the number 
of 
$\phi\in\param_{-\delta}$ for which $P_\phi+\lambda\subseteq \weylgon{\lambda}$.
(Part 1 of Theorem~\ref{prop:statement}).  
The multiplicity of $V_{\lambda+\mu-\delta}$ in $V_\lambda\otimes V_\mu$ equals the 
number of these that are also contained in 
$\weylgon{{-\mu}}+\lambda+\mu-\delta$
(Part 2 of Theorem~\ref{prop:statement}).  
Sharpness is because, 
given $\delta$, we can choose $\lambda$ and $\mu$ large enough that  
$\mm(\overline{\Gr_\lambda})\cap(\mm(\overline{\Gr_{-\mu}})+\lambda+\mu-\delta)=\mm(\overline{S_{\lambda-\delta}\cap T_\lambda})$ 
and that this contains only dominant weights.
\stopproof

\section{Hopf Algebra of MV-cycles}
\label{sec:hopf}

We close with a brief discussion of some related topics.  
This is mostly conjectural and is discussed in much more detail in~\cite{A}.  
There, we define a product on $\alg=\mbox{span}(\poly)$ and give 
a conjectural definition of a coproduct.  Conjecturally, $\alg$ 
is isomorphic 
the Hopf algebra of polynomial functions on the unipotent radical of a Borel 
subgroup of $\repG$.
Another loop Grassmannian approach to the (dual) Hopf algebra may be found in~\cite{FFKM}.

The product in $\alg$ is defined by a deformation of varieties over a curve, 
using an idea of Drinfeld's.  
There seem to be canonical generators and relations (described for some low 
rank groups in section~\ref{sec:examples}), but little is understood about them.  
Positive integer coefficients 
appear in the relations because they are multiplicities of 
irreducible components.  

The coproduct is not well understood.  
As an example, 
the coproduct and antipode for $\spfour$ are given in the following table.  
(It is only necessary to specify them on the generators since these functions 
are multiplicative.) 
\begin{center}
\vspace{5 mm}
\begin{tabular}{||c|c|l||}
\hline
Name & Antipode & \multicolumn{1}{c||}{Coproduct}  \\
$x$ & $S(x)$ & \multicolumn{1}{c||}{$\Delta(x)$}  \\
\hline
$1$ & $1$ & $1\otimes1$ \\
\hline
$a_1$ & $-a_1$ & $a_1\otimes1+1\otimes a_1$ \\
$a_2$ & $-a_2$ & $a_2\otimes1+1\otimes a_2$ \\
\hline
$b_1$ & $b_2$ & $b_1\otimes1+a_1\otimes a_2+1\otimes b_1$ \\
$b_2$ & $b_1$ & $b_2\otimes1+a_2\otimes a_1+1\otimes b_2$ \\
\hline
$c_1$ & $-c_2$ & $c_1\otimes1+2b_1\otimes a_2+a_1\otimes a_2^2+1\otimes c_1$ \\
$c_2$ & $-c_1$ & $c_2\otimes1+a_2^2\otimes a_1+2a_2\otimes b_2+1\otimes c_2$ \\
$c_3$ & $-c_3$ & $c_3\otimes1+b_2\otimes a_2 + a_2 \otimes b_1 + 1 \otimes c_3$ \\
\hline
$d_1$ & $d_1$ & $d_1\otimes1+c_1\otimes a_1+2b_1\otimes b_2+a_1\otimes c_2+1\otimes d_1$ \\
\hline
\end{tabular}
\end{center}
\vspace{\posttable}

Just as the product in $\alg$ corresponds to Minkowski sum of polytopes, the coproduct 
has a (conjectural) interpretation in terms of polytopes as well:  
Write $\Delta(x)=\sum k_{ij}x_i\otimes x_j^\prime$ 
where each $x_i$ is an MV-cycle relative to $\unip$ and each $x_j^\prime$ 
is an MV-cycle relative to $\unip^-$. 
If $k_{ij}\neq 0$ then (1)~$x_i$ and $x_j^\prime$ 
are contained in $x$; (2)~$x_i$ and $x_j^\prime$ are associated 
to the same weight in $x$, and this is their only point of intersection; 
(3)~$k_{ij}$ is a positive integer.  
The picture 
{\par 
\begin{center}
\includegraphics[width=12cm]{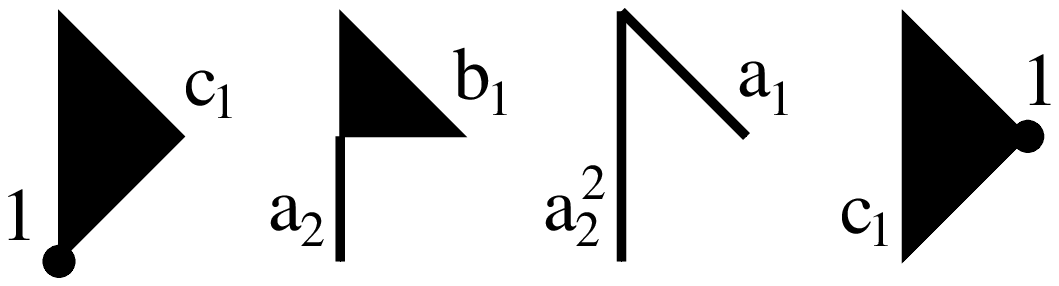} 
\end{center}
}
\noindent
illustrates this for
$\Delta(c_1)=c_1\otimes1+2b_1\otimes a_2+a_1\otimes a_2^2+1\otimes c_1$.

What are the meanings of the $k_{ij}$ and why are they positive?  
If $x$ is the largest MV-cycle in a stratum, 
they seem to give the intersection form in intersection homology~\cite{A}.  
If not (as in the above example) their meaning is mysterious.


\begin{center}
{\sc Acknowledgements}
\end{center}

I thank my thesis advisor, R.~MacPherson, and also D.~Nadler and I.~Mirkovi{\'c}, 
for many discussions 
on geometry in representation theory and for significant ideas that went into this research.  
I also thank J.~Conway, D.~Gaitsgory and K. Vilonen 
for discussions.  
Much of this work was done under an NSF fellowship.

\enlargethispage{1000pt}


\begin{thebibliography}{999999}

\bibitem[A]{A}
J. Anderson.
On Mirkovi{\'c} and Vilonen's Intersection Homology Cycles for the 
Loop Grassmannian.
PhD Thesis, Princeton University, 2000.

\bibitem[AM]{AM}
J. Anderson and I. Mirkovi{\'c}.  
Crystal Graphs via Polytopes
\textit{In preparation}, 2001.

\bibitem[AP]{AP}
M. Atiyah and A. Pressley.
Convexity and Loop Groups.  
\textit{Arithmetic and geometry, Vol.~II (Progr. Math \textbf{36})}
Birkh{\"a}user, Boston, 1983, pp.~33-63.

\bibitem[BBD]{BBD}
A. Beilinson, J. Bernstein, and P. Deligne. 
Faisceaux pervers.  Analyse et Topologie sur les Espaces Singuliers, 
volume 1, \textit{Ast{\'e}risque} \textbf{100}, 1982. 

\bibitem[BD]{BD}
A. Beilinson and V. Drinfeld.
Quantization of Hitchin's Integrable System and Hecke Eigensheaves. 
\textit{Work in progress}, 2000.


\bibitem[BG]{BG}
A. Braverman and D. Gaitsgory.  
Crystals via the affine Grassmannian. 
\textit{Duke Math. J.} \textbf{107} (2001) no.~3, pp.~561--575. 

\bibitem[B]{B}
M. Brion.
Sur l'image de l'application moment.  
\textit{Seminaire d'algebre Paul Dubreil et Marie-Paule Malliavin, Paris 1986 (Lecture notes in Mathematics \textbf{1296}).}
Springer, Berlin, 1987. pp.~177-192.

\bibitem[C]{C}
J. Carrell.
The Bruhat graph of a Coxeter group, a conjecture of Deodhar, and rational smoothness of Schubert varieties.
\textit{Proc. Symp. in Pure Math.} \textbf{56} (1994) pp.~53-61.


\bibitem[DM]{DM}
P. Deligne and J.Milne.  
Tannakian categories.
\textit{Hodge cycles, motives, and Shimura varieties 
(Lecture Notes in Mathematics \textbf{900}).}  
Springer 1982, \mbox{pp.~101-228}.

\bibitem[FFKM]{FFKM}
B. Feigin, M. Finkelberg, A. Kuznetsov, and I. Mirkovi{\'c}.
Loop Grassmannian Construction of $U(\check{n})$.
\textit{Work in progress}, 1999.


\bibitem[G]{VG}
V. Ginzburg.
Perverse Sheaves on a Loop Group and Langlands Duality.
\textit{Preprint}, 1990.

\bibitem[GM1]{GM1}
M. Goresky and R. MacPherson.
Intersection homology theory.
\textit{Topology} \textbf{19} (1980) no.~2, pp.~135-162.

\bibitem[GM2]{GM2}
M. Goresky and R. MacPherson.
Intersection homology II.
\textit{Invent. Math.} \textbf{72} (1983) no.~1, pp.~77-129.

\bibitem[GM3]{GM3}
M. Goresky and R. MacPherson.
On the topology of algebraic torus actions.  
\textit{Algebraic groups Utrecht 1986 (Lecture notes in Mathematics \textbf{1271}).}
Springer, Berlin-New York, 1987. pp.~73-90.




\bibitem[L]{L}
G. Lusztig. 
Singularities, character formulas, 
and a q-analog of weight multiplicities.
\textit{Ast{\'e}risque} \textbf{101-102} (1983), pp.~208-229.

\bibitem[M]{M}
R. MacPherson.
Intersection Homology and Perverse Sheaves.
\textit{Unpublished lecture notes}, 1991.

\bibitem[MV]{MV} 
I. Mirkovi{\'c} and K. Vilonen. 
Perverse Sheaves on affine Grassmannians and Langlands Duality.
\textit{Math. Res. Lett.} \textbf{7} (2000) no.~1, pp.~13-24.

\bibitem[PS]{PS}
A. Pressley and G. Segal. 
\textit{Loop Groups.}
Clarendon Press, Oxford, 1986.


\end{thebibliography}
\end{document}